\documentclass[11pt]{article}%
\usepackage{amsmath}
\usepackage{graphicx}
\usepackage{amsfonts}
\usepackage{amssymb}
\usepackage{a4wide}
\newtheorem{theorem}{Theorem}

\newtheorem{example}[theorem]{Example}

\begin{document}

\title{On some extension theorems for multi\-functions}
\author{C. Z\u{a}linescu\thanks{University Alexandru Ioan Cuza, Faculty of
Mathematics, Ia\c{s}i, Romania, email: \texttt{zalinesc@uaic.ro}.}}
\date{}
\maketitle

The extension of linear mappings dominated by convex multi\-functions is a
subject that interested quite many authors. When dealing with the continuous
case problems appear and several false result were published. In the
literature one can find correct versions of such results (of course, under
more stringent conditions). It is our aim here to point out some false results
published recently providing counterexamples. Note that in our paper
\cite{Zal:08b} we mentioned two wrong results. In this short paper we refer to
such results published after 2008.

We quote first \cite[Th.\ 2.2]{Du:09}.


\begin{quotation}
\noindent\textbf{Theorem 2.2} (Generalized Hahn Banach theorem in t.v.s.). Let
$X$ be a real t.v.s. and $C$ a convex subset of $X$. Let $(Y,K)$ be a real
order-complete t.v.s. with order cone $K$. Let $P:C\rightrightarrows Y$ be a
$K$-convex set-valued map with nonempty values. Let $X_{0}$ be a proper real
linear subspace of $X$ with $X_{0}\cap\operatorname*{cor}C\neq\emptyset$ and
$f:X_{0}\rightrightarrows Y$ be a continuous (u.s.c., l.s.c., respectively)
$K$-concave set-valued map with nonempty values satisfying

\medskip

$P(x)-f(x)\subset K$

\medskip

\noindent for all $x\in X_{0}\cap C$. Then there exists a set-valued map with
nonempty values $F:X\rightrightarrows Y$ such that

{(i)} $F$ is continuous (u.s.c., l.s.c., respectively);

{(ii)} $F$ is $K$-concave;

{(iii)} $F$ is an extension of $f$, i.e., $F(x)=f(x)$ for all $x\in X_{0}$;

{(iv)} $P(x)-F(x)\subset K$ for all $x\in C$.
\end{quotation}


The fact that $P$ is $K$-convex means that $\left\{  (x,y+k)\mid x\in X,\ y\in
P(x),\ k\in K\right\}  $ is convex; $f$ is $K$-concave means that $f$ is
$(-K)$-convex.

\begin{example}
\label{ce-dut22}Take $X$ an infinite dimensional real normed vector space,
$Y:=\mathbb{R}$ endowed with the usual topology and $K:=\mathbb{R}%
_{+}:=[0,\infty).$ Of course $(Y,K)$ is a real order-complete t.v.s. Consider
$\varphi:X\rightarrow\mathbb{R}$ a linear non continuous function and
$P:X\rightrightarrows\mathbb{R}$ defined by $P(x):=\{\varphi(x)\}$. Consider
also $X_{0}:=\{0\}\subset X$ and $f:X_{0}\rightrightarrows\mathbb{R}$ defined
by $f(0):=\{0\}.$ With these data the hypothesis of \cite[Th.\ 2.2]{Du:09}
holds but its conclusion is not valid.
\end{example}

Clearly $P$ is $K$-convex, and $f$ is continuous and $K$-concave. Moreover,
$C:=X=\operatorname*{cor}C,$ whence $0\in X_{0}\cap\operatorname*{cor}C$, and
$P(x)-f(x)=\{0\}\subset K$ for $x\in X_{0}\cap C=\{0\}.$ Applying the theorem
above, that is \cite[Th.\ 2.2]{Du:09}, we find $F:X\rightrightarrows
\mathbb{R}$ which is $K$-concave, continuous, $F(0)=\{0\}$ and $\emptyset
\neq\{\varphi(x)\}-F(x)\subset K,$ that is $F(x)\subset\{\varphi
(x)\}-\mathbb{R}_{+},$ for all $x\in X.$ Since $\{0\}=F(0)\subset(-1,1)$ and
$F$ is upper semi\-continuous at $0,$ there exists $r>0$ such that
$\emptyset\neq F(x)\subset(-1,1)$ for every $x\in U:=\{u\in X\mid\left\Vert
u\right\Vert \leq r\}.$ Taking $x\in U$ and some $y\in F(x)$ we obtain that
$-1\leq y\in\{\varphi(x)\}-\mathbb{R}_{+};$ hence $\varphi(x)\geq-1$ for every
$x\in U.$ Therefore, $\varphi$ is continuous, a contradiction.

\medskip

The above example is also a counterexample for \cite[Cor.\ 2.2]{Du:09}.
Indeed, take the same $P$ and $u=0\in X,$ $v=0\in\mathbb{R}.$

Let us quote now \cite[Th.\ 3.2]{Du:09}.

\begin{quotation}
\noindent\textbf{Theorem 3.2} Let $X$ be a real t.v.s. and $C$ a convex subset
of $X$. Let $(Y,K)$ be a real order-complete t.v.s. with order cone $K$. Let
$f:C\rightrightarrows Y$ be a $K$-convex set-valued map with nonempty values.
Let $X_{0}$ be a proper real linear subspace of $X$ with $X_{0}\cap
\operatorname*{cor}C\neq\emptyset$ and $P:X_{0}\rightrightarrows Y$ be a
continuous (u.s.c., l.s.c., respectively) $K$-concave set-valued map with
nonempty values satisfying

\medskip

$P(x)-f(x)\subset K$

\medskip

\noindent for all $x\in X_{0}\cap C$. Then there exists a set-valued map with
nonempty values $G:X\rightrightarrows Y$ such that

{(i)} $G$ is continuous (u.s.c., l.s.c., respectively);

{(ii)} $G$ is $K$-concave;

{(iii)} $G$ is an extension of $P$, i.e., $G(x)=P(x)$ for all $x\in X_{0}$;

{(iv)} $G(x)-f(x)\subset K$ for all $x\in C$.
\end{quotation}

\begin{example}
\label{ce-dut32}Take $X,$ $X_{0},$ $(Y,K)$ and $\varphi$ as in Example
\ref{ce-dut22}. Let $f:X\rightrightarrows\mathbb{R}$ be defined by
$f(x):=\{\varphi(x)\}$, and $P:X_{0}\rightrightarrows\mathbb{R}$ be defined by
$P(0):=\{0\}.$ With these data the hypothesis of \cite[Th.\ 3.2]{Du:09} holds
but its conclusion is not valid.
\end{example}

Clearly, all the hypotheses of \cite[Th.\ 3.2]{Du:09} hold. Applying this
theorem we get $G:X\rightrightarrows\mathbb{R}$ with $G(0)=\{0\},$
$\emptyset\neq G(x)\subset\{\varphi(x)\}-\mathbb{R}_{+}$ for all $x\in X$ and
$G$ continuous. Since $G(0)=\{0\}\subset(-1,1)$ and $G$ is upper
semi\-continuous at $0,$ there exists $r>0$ such that $G(x)\subset(-1,1)$ for
every $x\in U:=\{u\in X\mid\left\Vert u\right\Vert \leq r\}.$ Taking $x\in U$
and some $y\in G(x)$ we obtain that $1\geq y\in\{\varphi(x)\}+\mathbb{R}_{+};$
hence $\varphi(x)\leq1$ for every $x\in U.$ Therefore, $\varphi$ is
continuous, a contradiction.

\medskip

The above example is also a counterexample for \cite[Cor.\ 3.2]{Du:09}.
Indeed, take the same $f$ and $u=0\in X,$ $v=0\in\mathbb{R}.$

\medskip

In \cite{LiGuo:09} one says: ``From Theorem 13 in [11] and its proof, we have
the following result''. (Surely, it is $T_{0}\in L(X_{0},Y)$ instead of
$T_{0}\in L(X,Y)$.)

Let us quote \cite[Th.\ 2.1]{LiGuo:09}.

\begin{quotation}
\noindent\textbf{Theorem 2.1} ([11]). Let $(Y,K)$ have the least upper bound
property, $F:X\rightarrow2^{Y}$ be a $K$-convex set-valued mapping,
$X_{0}\subset X$ be a linear subspace and $T_{0}\in L(X,Y)$. Suppose that
$0\in\operatorname*{int}(\operatorname*{dom}F-X_{0})$ and $T_{0}(x)\not >  y$
for all $(x,y)\in\operatorname*{gph}F\cap(X_{0}\times Y)$. Then there exists
$T\in L(X,Y)$ such that $T|_{X_{0}}=T_{0}$ and $T(x)\not >  y$ for all
$(x,y)\in\operatorname*{gph}F$. Moreover, $T$ can be defined by $T(x)=T_{0}%
(x_{0})-\left\langle x_{1},x^{\ast}\right\rangle \overline{y}$, where
$\overline{y}\in\operatorname*{int}K$, $x^{\ast}\in X^{\prime}$ and the linear
subspace $X_{1}\subset X$ with $X=X_{0}\oplus X_{1}$ are fixed, and
$x=x_{0}+x_{1}$ with $x_{0}\in X_{0}$ and $x_{1}\in X_{1}$.
\end{quotation}


In the theorem above $X,Y$ are real linear topological spaces, and $L(X,Y)$ is
the set of all linear continuous operators from $X$ into $Y$; furthermore,
$X^{\prime}=L(X,\mathbb{R})$, and $K\subset Y$ is a proper pointed closed
convex cone with nonempty interior $\operatorname*{int}K$.

\begin{example}
\label{ce-LG21}Take $X,$ $X_{0},$ $(Y,K)$ and $\varphi$ as in
Example \ref{ce-dut22}. Let $F:X\rightrightarrows\mathbb{R}$ be
defined by $F(x):=\{\varphi(x)\}$ and
$T_{0}:X_{0}\rightarrow\mathbb{R}$ be defined by $T_{0}(0):=\{0\}.$
With these data the hypothesis of \cite[Th.\ 2.1]{LiGuo:09} holds
but its conclusion is not valid.
\end{example}

Clearly, all the hypotheses of \cite[Th.\ 2.1]{LiGuo:09} are
satisfied; applying it we get $T\in L(X,\mathbb{R})=X^{\prime}$ such
that $T(x)\not > y$ for all $(x,y)\in\operatorname*{gph}F$, that is
$T(x)\leq\varphi(x)$ for all $x\in X.$ It follows that $T=\varphi,$
and so we get the contradiction that $\varphi$ is continuous.

\medskip

Example \ref{ce-LG21} is also a counterexample for \cite[Th.\
3.1]{LiGuo:09} quoted below; just take $(x_{0},y_{0}):=(0,0).$


\begin{quotation}
\noindent\textbf{Theorem 3.1}. Let $(Y,K)$ have the least upper bound
property, $F:X\rightarrow2^{Y}$ be a $K$-convex set-valued mapping with
$\operatorname*{int}(\operatorname*{dom}F)\neq\emptyset$. If there exist
$x_{0}\in\operatorname*{int}(\operatorname*{dom}F)$ and $y_{0}\in F(x_{0})$
such that $F(x_{0})\cap(y_{0}-\operatorname*{int}K)=\emptyset$, then
$\partial^{Y-W}F(x_{0},y_{0})\neq\emptyset$. Moreover, there is a Y-weak
subgradient $T$ of $F$ at $(x_{0},y_{0})$ satisfying that for every
$x\in\operatorname*{dom}F$,

$T(x-x_{0})\notin-\operatorname*{int}K\Leftrightarrow T(x-x_{0})\in K.$
\end{quotation}


Here, for $(x_{0},y_{0})\in\operatorname*{gph}F,$
\[
\partial^{Y-W}F(x_{0},y_{0}):=\{T\in L(X,Y)\mid y-y_{0}-T(x-x_{0}%
)\notin-\operatorname*{int}K\ \forall(x,y)\in\operatorname*{gph}F\}.
\]

\medskip

Example \ref{ce-LG21} is a counterexample for \cite[Th.\ 3.1]{GuoLiTeo:12}
quoted below, too; just take $(\overline{x},\overline{y}):=(0,0).$


\begin{quotation}
\noindent\textbf{Theorem 3.1} Let $(Y,K)$ have the order-complete
property, $F:X\rightarrow2^{Y}$ be a $K$-convex set-valued mapping
with $\operatorname*{int}(\operatorname*{dom}F)\neq\emptyset$. If
there exist
$\overline{x}\in\operatorname*{int}(\operatorname*{dom}F)$ and
$\overline {y}\in F(\overline{x})$ such that
$F(\overline{x})-\overline{y}\subset K$, then $\partial
F(\overline{x},\overline{y})\neq\emptyset$.
\end{quotation}


Here, $X,Y$ are real normed spaces, $L(X,Y)$ is the set of all linear
continuous operators from $X$ to $Y$, $K\subset Y$ is a proper pointed closed
convex cone with nonempty interior $\operatorname*{int}K,$ and for
$(\overline{x},\overline{y})\in\operatorname*{gph}F,$
\[
\partial F(\overline{x},\overline{y}):=\{T\in L(X,Y)\mid y-\overline
{y}-T(x-\overline{x})\in K\ \forall(x,y)\in\operatorname*{gph}F\}.
\]

\end{document}